\documentclass{amsart}
\usepackage{amssymb}
\usepackage{graphicx}

\newcommand{\face}[1]{\xrightarrow{#1}}
\newcommand{\inv}{^{-1}}
\newcommand{\co}{\colon\thinspace}
\newcommand{\inter}{\textrm{int}}

\newtheorem{thm}{Theorem}

\newtheorem{cor}[thm]{Corollary}

\newtheorem{ex}[thm]{Example}
\newtheorem{lemma}[thm]{Lemma}

\newtheorem{remark}[thm]{Remark}

\newcommand{\sections}{\renewcommand{\thethm}{\thesection.\arabic{thm}}
           \setcounter{thm}{0}}

\newcommand{\subsections}{\setcounter{thm}{0}\renewcommand{\thethm}
           {\thesubsection.\arabic{thm}}}
\newcommand{\nosubsections}{\renewcommand{\thethm}
{\thesection.\arabic{thm}}
           \setcounter{thm}{0}}

\newcommand{\nosubsubsections}{\setcounter{thm}{0}
           \renewcommand{\thethm}{\thesubsection.\arabic{thm}}}

\newcommand{\bR}{\mathbb{R}}
\newcommand{\bQ}{\mathbb{Q}}
\newcommand{\bZ}{\mathbb{Z}}
\newcommand{\Zpm}{\mathbb{Z}\setminus\{0\}}

\newcommand{\za}{\alpha}
\newcommand{\zb}{\beta}

\newcommand{\zd}{\delta}
\newcommand{\ze}{\epsilon}

\newcommand{\zs}{\sigma}
\newcommand{\zt}{\tau}
\newcommand{\zv}{\varphi}

\newcommand{\cD}{\mathcal{D}}
\newcommand{\cN}{\mathcal{N}}
\newcommand{\cV}{\mathcal{V}}

\newcommand{\link}{\textrm{link}}
\newcommand{\mul}{\textrm{mul}}

\begin{document}
\title{Bitwist 3-manifolds}

\author{J. W. Cannon}
\address{Department of Mathematics\\ Brigham Young University\\ Provo, UT
84602\\ U.S.A.} \email{cannon@math.byu.edu}

\author{W. J. Floyd}
\address{Department of Mathematics\\ Virginia Tech\\
Blacksburg, VA 24061\\ U.S.A.} \email{floyd@math.vt.edu}
\urladdr{http://www.math.vt.edu/people/floyd}

\author{W. R. Parry}
\address{Department of Mathematics\\ Eastern Michigan University\\
Ypsilanti, MI 48197\\ U.S.A.} \email{walter.parry@emich.edu}

\date{\today}

\begin{abstract} Our earlier twisted-face-pairing construction showed
how to modify an arbitrary orientation-reversing face-pairing $\ze$ on
a faceted 3-ball in a mechanical way so that the quotient is
automatically a closed, orientable $3$-manifold.  The modifications
were, in fact, parametrized by a finite set of positive integers,
arbitrarily chosen, one integer for each edge class of the original
face-pairing. This allowed us to find very simple face-pairing
descriptions of many, though presumably not all, 3-manifolds.

Here we show how to modify the construction to allow negative
parameters, as well as positive parameters, in the
twisted-face-pairing construction. We call the modified construction
the bitwist construction. We prove that all closed connected
orientable 3-manifolds are bitwist manifolds.  As with the twist
construction, we analyze and describe the Heegaard splitting naturally
associated with a bitwist description of a manifold.
\end{abstract}

\maketitle

\sections

\section{Introduction}\label{sec:intro}\nosubsections

In a series (\cite{introtfp,twist,heegaard}) of papers, we described
and analyzed a simple construction of 3-manifolds from face-pairings.
If $\ze$ is an orientation-reversing edge-pairing on a polygonal disk
$D$, then the quotient space $D/\ze$ is always a surface. But if $\ze$
is an orientation-reversing face-pairing on a faceted 3-ball $P$, the
quotient $P/\ze$ is not generally a 3-manifold.  (See, for example,
Section 2.7 of \cite{DT}.) For the twist construction one chooses a
positive integer, called the multiplier, for each edge cycle
(equivalence class of an edge under the action of $\ze$). By
subdividing each edge into the product of its multiplier and the size
of its edge cycle and then precomposing $\ze$ with a twist, one
obtains a new faceted 3-ball $Q$ and orientation-reversing face
pairing $\zd$. The fundamental result of the constuction is that
$Q/\zd$ is always a 3-manifold. Papers \cite{introtfp} and
\cite{twist} give the basic details of the construction. The
construction is analyzed further in \cite{heegaard}, and Heegaard
diagrams and surgery diagrams are given for twisted face-pairing
manifolds.

In this paper we give a modified construction which we call the
bitwist construction. The basic setup is the same, but we allow the
edge cycle multipliers to be positive or negative.  Allowing twisting
in different directions leads to problems in defining the new
face-pairing $\zd$, but one can resolve this by the appropriate
insertion of ``stickers'' in the faces of the new faceted
3-ball $Q$. In Section~\ref{sec:prelimex}
we give a simple preliminary example to show how stickers are used in
the construction. Following this, we give the general construction in
Section~\ref{sec:const}.  As with the twist construction, the
3-manifolds constructed from the bitwist construction naturally have a
cell structure with a single vertex. One can easily give presentations
for fundamental groups of bitwist manifolds as in \cite[Section~4]{twist},
but the homology results of \cite[Secton~6]{twist} do not generally
hold for bitwist manifolds. Since $Q/\zd$ has a single
vertex, some of the
results from the twist construction apply directly to the bitwist
construction. In particular, the construction of Heegaard diagrams and
framed surgery descriptions from \cite{heegaard} are valid for the
bitwist construction.  This is developed in
Section~\ref{sec:heegaard}. If $L$ is a corridor complex link for an
orientation-reversing face pairing $\ze$ on a faceted 3-ball $P$ and
$\mul$ is a multiplier function for $(P,\ze)$, then the bitwist
manifold $M(P,\ze,\mul)$ is obtained by framed surgery on $L$, where
the face components have framing $0$ and an edge component has framing
the sum of its blackboard framing and the reciprocal of the multiplier
of its edge cycle.

After making the leap to negative multipliers, it is natural to
inquire about multipliers with value 0.  Allowing edge cycle
multipliers to be 0 amounts to collapsing every edge with multiplier 0
to a point and applying the construction to the resulting complex.  In
terms of our surgery description, this amounts to deleting from our
framed link every component with framing $\infty$, an operation which
does not change the resulting manifold.  Collapsing edges in general
leads to complexes which are no longer 3-balls -- they are cactoids.
While we actually do find face-pairings on cactoids interesting and we
do temporarily allow multipliers to be 0 in the proof of
Theorem~\ref{thm:lens}, for the present we content ourselves with
nonzero multipliers.

The framed surgery descriptions are a primary motivation for
developing the bitwist construction. In order to realize 3-manifolds
as twisted face-pairing manifolds or bitwist manifolds, one wants to
be able to change the framings of the edge components. Suppose $L$ is
a corridor complex link for a twisted face-pairing manifold. We still
get a twisted face-pairing manifold if we replace the framing of each
edge component by its blackboard framing plus an arbitrary positive
rational number.  In Section~\ref{sec:framings} we show that using the
bitwist construction, we still get a bitwist manifold if we replace
the framing of each edge component by its blackboard framing plus an
arbitrary rational number. This ability to change the signs of the
rational numbers gives extra power to the construction.  Using this,
we show in Section~\ref{sec:allarebitwist} that every closed connected
orientable 3-manifold is a bitwist manifold.

\section{A preliminary example}\label{sec:prelimex}\nosubsections

\begin{figure}
\centerline{\includegraphics{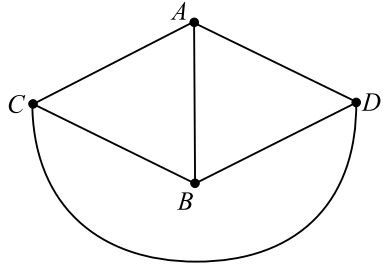}} \caption{The tetrahedron
$P$.} \label{fig:teta}
\end{figure}

We give a preliminary example to indicate the construction. We start
with a simple model face-pairing $\ze$ that was considered in
Section~2 of \cite{introtfp}. Our faceted 3-ball $P$  is a
tetrahedron with vertices $A$, $B$, $C$, and $D$, as shown in
Figure~\ref{fig:teta}. We consider $P$ as an oriented 3-ball, and
for convenience give it an orientation so that in the induced
orientation on the boundary of $P$ the boundary of each 2-cell is
oriented clockwise.

The model face pairing $\ze$ identifies the triangles $ABC$ and $ABD$
by reflection in the common edge $AB$, and it identifies $ACD$ and
$BCD$ by reflection in the common edge $CD$. In the permutation
notation of \cite{introtfp}, $\ze$ is given as follows:
\begin{equation*}
\ze_1:\begin{pmatrix} A   &B  &C  \\
               A   &B  &D  \end{pmatrix}\quad
\ze_2:\begin{pmatrix} A   &C  &D  \\
               B   &C  &D \end{pmatrix} .
\end{equation*}

There are three edge cycles, as follows:
\begin{equation*}
AB\face{\ze_1}AB
\end{equation*}
  \begin{equation*}
BC\face{\ze_1}BD\face{\ze_2\inv}AD\face{\ze_1\inv}AC\face{\ze_2}BC
\end{equation*}
  \begin{equation*}
CD\face{\ze_2}CD.
\end{equation*}
The first edge cycle $[AB]$ has length $\ell([AB])=1$, the second edge
cycle $[BC]$ has length $\ell([BC])=4$, and the third edge cycle
$[CD]$ has length $\ell([CD])=1$. In the twisted face-pairing
construction, for each edge cycle one chooses a positive integer
$\mul([e])$ called the \textbf{multiplier}. For the bitwist construction,
one chooses a nonzero integer $\mul([e])$, still called the
multiplier, for each edge cycle. We use the cycle lengths and the
absolute values of the multipliers to determine how to subdivide the
edges of $P$.  The sign of the multiplier indicates the direction in
which we twist edges in the edge cycle $[e]$. If all of the multipliers
have the same sign, then we have the twist construction. For this
example, we choose $\mul([AB])=-1$, $\mul([BC])=1$, and $\mul([CD])=1$.

\begin{figure}
\centerline{\includegraphics{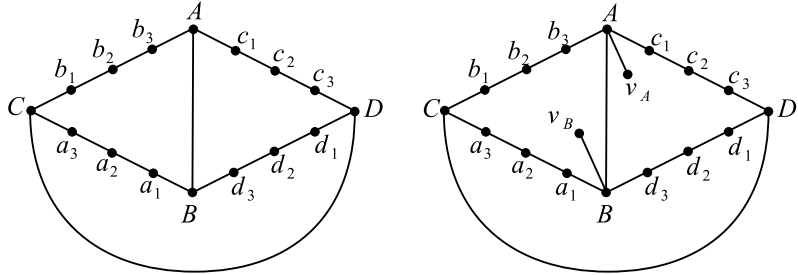}} \caption{The subdivisions
$Q'$ and $Q$ of $P$.} \label{fig:tetasub}
\end{figure}

We are now ready to replace $P$ by its subdivision $Q$. We subdivide
every edge $e$ of $P$ into $\ell([e]) \cdot |\mul([e])|$ subedges. We
perform these subdivisions so that the face-pairing $\ze$ takes
subedges to subedges. Let $Q'$ be the resulting faceted $3$-ball. We
need to perform a further modification if the multipliers do not all
have the same sign.  Let $f$ be a face of $P$.  Suppose $v$ is a
vertex of $P$ in $f$. Let $e_1$ be the edge of $P$ in $f$ with
terminal vertex $v$ and let $e_2$ be the edge of $P$ in $f$ with
initial vertex $v$. If $\mul([e_1]) < 0$ and $\mul([e_2]) > 0$, then
we add a \textbf{sticker} (think straight pin with spherical head) to
$f$ at $v$. That is, we add a new vertex in the interior of $f$ and
join it to $v$ by an edge in $f$. The faceted 3-ball obtained from $P$
by adding stickers to $Q'$ as described above is the subdivision
$Q$. Figure~\ref{fig:tetasub} shows the subdivisions $Q'$ and $Q$ for
this example.

We define a bitwisted face-pairing $\zd$ on $Q$ as follows:
  \begin{equation*}
\zd_1:\left(\begin{array}{ccccccccccc}
b_3&A&B&v_B&B&a_1&a_2&a_3&C&b_1&b_2\\
            A&v_A&A&B&d_3&d_2&d_1&D&c_3&c_2&c_1\end{array}\right);
\end{equation*}
  \begin{equation*}
\zd_2:\left(\begin{array}{ccccccccc} c_1&A&b_3&b_2&b_1&C&D&c_3&c_2\\
           B&a_1&a_2&a_3&C&D&d_1&d_2&d_3\end{array}\right).
\end{equation*}
The underlying idea is that we precompose $\ze$ with a twist in the
positive direction on an edge which is a subedge of an original edge
with positive multiplier, and we precompose $\ze$ with a twist in the
negative direction on an edge which is a subedge of an original edge
with negative multiplier. This is not well defined on $Q'$ since
adjacent original edges can have multipliers of different signs, but one
can make it well defined on $Q$.

\begin{figure}
\centerline{\includegraphics{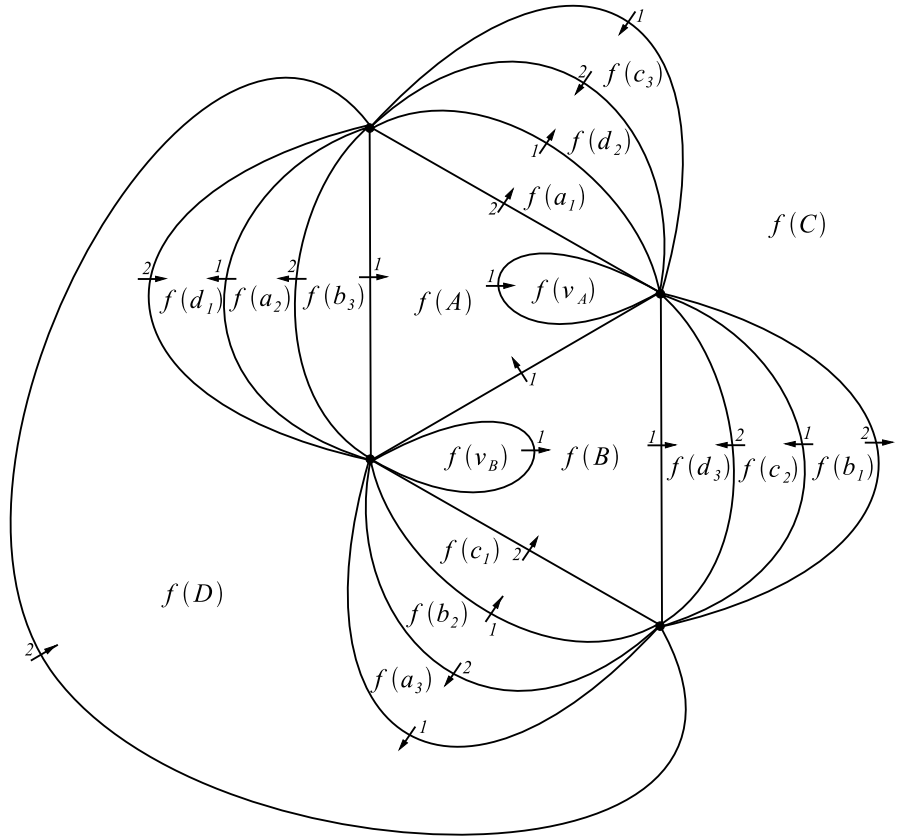}} \caption{The link of the
vertex of $M$.} \label{fig:tetalink}
\end{figure}

\begin{figure}
\centerline{\includegraphics{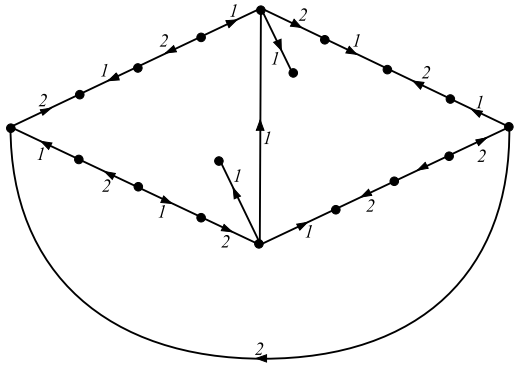}} \caption{The dual
subdivision.} \label{fig:tetadual}
\end{figure}

Let $M = Q/\sim$, where $\sim$ is the equivalence relation on $Q$
generated by the face-pairing $\zd$. The computation below shows
that $M$ has two 1-cells and a single 0-cell.
\begin{equation*}
\begin{split} b_3A\face{\zd_1}Av_A\face{\zd_1\inv} & BA
\face{\zd_1\inv}v_BB \face{\zd_1} Bd_3\face{\zd_2\inv}c_1c_2
  \face{\zd_1\inv} b_2b_1\face{\zd_2} \\  & a_3C\face{\zd_1}Dc_3
  \face{\zd_2}d_1d_2\face{\zd_1\inv}a_2a_1\face{\zd_2\inv}b_3A \\
  b_1C\face{\zd_2}CD\face{\zd_2} & Dd_1\face{\zd_1\inv}a_3a_2
  \face{\zd_2\inv}b_2b_3\face{\zd_1} \\ & c_1A\face{\zd_2}Ba_1
  \face{\zd_1}d_3d_2\face{\zd_2\inv}c_2c_3\face{\zd_1\inv}b_1C
\end{split}
\end{equation*}
Since $M$ has two 2-cells and a single 3-cell, $\chi(M) = 0$ and so
$M$ is a 3-manifold. Figure~\ref{fig:tetalink} shows the link of the
vertex of $M$. As for the twist construction, $M$ can also be obtained
as the quotient under the face pairings of a dual faceted 3-ball
$\partial Q^{*}$, and the boundary of $\partial Q^{*}$ is cellularly
isomorphic to the dual of the link shown in Figure~\ref{fig:tetalink}.
The subdivision of $\partial Q^{*}$ is shown in
Figure~\ref{fig:tetadual}. It is easy to see from Figure~\ref{fig:tetadual}
or from the display above that
\begin{eqnarray*}
\pi_1(M)& \cong & \langle x,y\co x x^{-1} x^{-1} x y^{-1}x^{-1}y x
y x^{-1} y^{-1}, y y x^{-1} y^{-1} x y x y^{-1} x^{-1} \rangle\\
& \cong & \langle x,y\co y^{-2}x^{-1}y x y x^{-1},
y^2 x^{-1} y^{-1} x y x y^{-1} x^{-1} \rangle.\\
\end{eqnarray*}

\section{The bitwist construction}\label{sec:const}\nosubsections

We now give the main construction. In \cite{twist} we defined a
faceted 3-ball to be a regular CW complex. Here we follow the  more
general definition of a faceted 3-ball $P$ given in \cite{heegaard}.
In particular, we do not assume that the 2-cells in $\partial P$ are
regular. As in \cite{heegaard}, a \textbf{faceted 3-ball} $P$ is an
oriented CW complex such that $P$ is a closed 3-ball, there is a
single 3-cell and its interior is $\inter(P)$, and $\partial P$ does
not consist solely of a 0-cell and a 2-cell. It follows from this
that for each 2-cell $f$ of $P$, there is a CW structure on a closed
disk $F_f$ such that i) $F_f$ has a single 2-cell and its interior is
$\inter(F_f)$ and ii) there is a continuous cellular map $\varphi\co
F_f\to f$ whose restriction to each open cell is a homeomorphism.

Still following \cite{heegaard}, given a faceted 3-ball $P$ we
construct a subdivision $P_s$ of $P$ by barycentrically subdividing
$\partial P$. The faceted 3-ball $P_s$ is a regular CW complex and
each 2-cell of $P_s$ is a triangle. Since the 2-cells of $P$ may not
be regular, a face pairing $\ze$ on $P$ is technically a matching of
the faces of $P$ together with a face pairing on $P_s$ which is
compatible with it. We still denote by $\ze$ the face pairing on
$P_s$. We assume as before that our face-pairings reverse orientation
and satisfy the face-pairing compatibility condition.

Suppose $P$ is a faceted 3-ball and $\ze$ is a face-pairing on $P$.
We refer to $\ze$ as a \textbf{model face-pairing}.  There is an
equivalence relation $\sim$ defined on the edges of $P$ that is
generated by the relation $e_1 \sim e_2$ if $e_2$ is the image of
$e_1$ under some element of $\ze$; the equivalence classes are called
\textbf{edge cycles}. If $[e]$ is an edge cycle, we denote its
cardinality by $\ell([e])$ and call it the \textbf{length} of $[e]$. In
addition to $(P,\ze)$, the input for the bitwist construction consists
of a \textbf{multiplier function}. The multiplier function is a function
$\textrm{mul}\co\{\textrm{edge cycles}\}\to \Zpm$.  An edge $e$ is
\textbf{positive} if $\mul([e])>0$ and is \textbf{negative} if
$\mul([e])<0$.

Suppose we are given a face-pairing $(P,\ze)$ together with a
multiplier function $\mul$. We create a subdivision $Q$ of $P$ in
two stages. The first stage consists of subdividing each edge $e$ of
$P$ into $\ell([e])\cdot |\mul([e])|$ subedges to get a subdivision
$Q'$ of $P$, and forming the subdivision $Q_s'$ of $Q'$ by
barycentrically subdividing $\partial Q'$. We perform these subdivisions so
that $\ze$ defines a face-pairing $\ze'$ on $Q'_s$.
The second stage of our construction of $Q$ consists of adding
stickers at some of the corners of the faces of $Q'$. Suppose $f$ is a
face of $P$, and consider a corner of $f$ at a vertex $v$ with edges
$e$ and $e'$, labeled such that $e'$ precedes $e$.  Suppose that $e'$
is a negative edge and $e$ is a positive edge.  Let $a\subseteq f$ be
the edge of $Q'_s$ which bisects this corner.  To $Q'$ we add a
barycenter $u$ of $a$ and the subedge of $a$ joining $u$ and $v$.
This subedge of $a$ is a sticker.  We continue with this process for
all of the corners of all of the faces of $P$. The result is a faceted
3-ball $Q$ which is obtained from $P$ by subdividing edges and adding
stickers.

\begin{figure}[!ht]
\centerline{\includegraphics{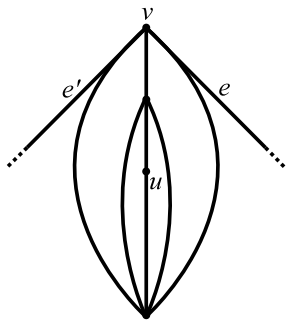}} \caption{The subdivision of a
face of $Q_s$ near a sticker.}
\label{fig:split}
  \end{figure}

As for $P$ and $Q'$, we form the subdivision $Q_s$ from $Q$ by
barycentrically subdividing $\partial Q$. We do this so that $Q_s$ is
a subdivision of $Q'_s$.  If $f$ is a face of $P$, we will still use
the name $f$ for the corresponding face in $Q$; to cut down on the
confusion, we will refer to edges of $P$ in $f$ as \textbf{original
edges} and to vertices of $P$ in $f$ as \textbf{original vertices}. Note
that $Q_s$ can be obtained from $Q_s'$ by splitting certain edges
which connect original vertices to barycenters of faces and then for
each split edge inserting a digon decomposed into four triangles. See
Figure~\ref{fig:split}, where the edge of $Q$ joining $u$ and $v$ is a
sticker.  In particular, there is a correspondence between faces of
$Q_s'$ and faces of $Q_s$ that do not contain subedges of stickers.

We next define a bitwisted face-pairing $\zd$ on $Q_s$. The
orientation on $P$, and hence on $Q$ and $Q_s$, determines a cyclic
order on the boundary of each face $f$ of $Q$ and hence a cyclic
order on the faces of the subdivision $f_s$.

\begin{figure}[!ht]
\centerline{\includegraphics{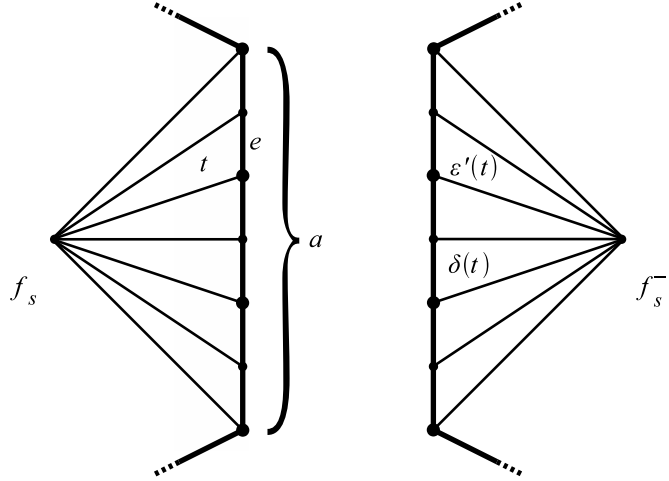}} \caption{Defining the
bitwisted face-pairing $\zd$.}
\label{fig:deltat1}
  \end{figure}

\begin{figure}[!ht]
\centerline{\includegraphics{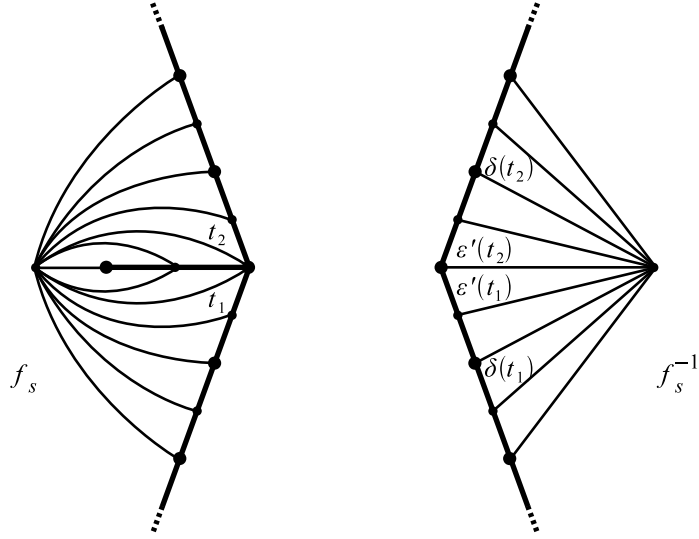}} \caption{Defining
$\zd$ near a sticker.}
\label{fig:deltat2}
  \end{figure}

Let $f$ be a face of $Q$, and let $e$ be an edge of $f_s$ which is
part of an original edge $a$ of $P$. See Figure~\ref{fig:deltat1},
which shows part of $f_s$ and $f_s^{-1}$ for some face $f$ of $Q$ with
positive original edge $a$.  The vertices and edges of $f$ and
$f^{-1}$ are drawn thick for emphasis.  Let $t$ be a face of the
subdivision $f_s$ which contains $e$. If $a$ is a positive edge, let
$\zd(t)$ be the face of $f^{-1}_s$ which is the second face before the
face $\ze'(t)$ of $f^{-1}_s$. If $a$ is a negative edge, let $\zd(t)$
be the face of $f^{-1}_s$ which is the second face after the face
$\ze'(t)$ of $f^{-1}_s$.  Figure~\ref{fig:deltat2} shows $\zd(t_1)$
and $\zd(t_2)$ for certain faces $t_1$ and $t_2$ of $f_s$ for the case
in which $f$ has a sticker.  The faces $t_1$ and $t_2$ both contain an
original vertex which is contained in the sticker.  Note that in
$f_s^{-1}$ from $\zd(t_1)$ to $\zd(t_2)$ in the positive direction
there are four faces corresponding to the four faces of $f_s$ which
contain a subedge of the sticker.  It follows that the definition of
$\zd$ can be extended to a face-pairing between $f_s$ and
$f^{-1}_s$. Doing this for each face defines a face-pairing $\zd$ on
$Q$. Unless the sign of $\mul$ is constant, this will not define a
face-pairing on $Q'_s$. In effect we are using the sign of $\mul$ to
determine which direction to twist each face of $Q_s$; the stickers
enable us to make this well defined.

We denote by $M(P,\ze,\mul)$ the quotient space of $Q$ under the
equivalence relation generated by $\zd$.

\begin{thm}\label{thm:main1} Let $P$ be a faceted 3-ball, let $\ze$ be
an orientation-reversing face-pairing on $P$ and let $\mul$ be a
multiplier function for $(P,\ze)$.  Then $M = M(P,\ze,\mul)$ is a
closed 3-manifold.  Furthermore, as a cell complex $M$ has just one
vertex.
\end{thm}
\begin{proof} The proof of the first assertion is an
Euler-characteristic argument analogous to the argument in
\cite{introtfp}. To prove that $M$ is a closed 3-manifold, it
suffices to show that $\chi(M)=0$.  We do this by determining the
number of cells in $M$ of every dimension.  It is clear that $M$ has
one 3-cell and that the number of 2-cells is the number of pairs of
faces of $Q$.  So to prove Theorem~\ref{thm:main1}, it suffices to
prove that $M$ has one 0-cell and that the number of 1-cells is the
number of pairs of faces of $Q$.

Every edge of $Q$ is either a sticker or a subedge of an original
edge.  The discussion involving Figure~\ref{fig:deltat2} shows that
the image under $\zd$ of a sticker contained in a face $f$ of $Q$
consists of two edges of $f^{-1}$.  One of these edges of $f^{-1}$ is
a terminal subedge of a positive original edge and one is an initial
subedge of a negative original edge.  The discussion involving
Figure~\ref{fig:deltat1} implies that every edge of $Q$ contained in
an original edge is equivalent to an edge $e$ of a face $f$ of $Q$
such that either $e$ is the terminal subedge of a positive original
edge of $f$ or $e$ is the initial subedge of a negative original edge
of $f$.  We conclude that every edge of $Q$ is equivalent to an edge
$e$ of a face $f$ such that either $e$ is the terminal subedge of a
positive original edge of $f$ or $e$ is the initial subedge of a
negative original edge of $f$.  Also, if $f$ is a face of $Q$ with a
positive original edge $e$ followed immediately by a negative original
edge $e'$, then the terminal subedge of $e$ is equivalent to the
initial subedge of $e'$ by means of a sticker.  Moreover every vertex
of $Q$ is equivalent to an original vertex.

Now let $e_0$ be an edge of a face $f_0$ of $Q$ such that $e_0$ is the
terminal subedge of a positive original edge of $f_0$.  Also suppose
that the original edge of $f_0$ immediately following $e_0$ is
positive.  By considering the $\zd$-orbit of $e_0$ we obtain edges
$e_1,\dotsc,e_n$ of faces $f_1,\dotsc,f_n$ of $Q$ and original edges
$e'_1,\dotsc,e'_n$ with the following properties.
  \begin{equation*}
\begin{gathered}
n=\ell([e'_i])\left|\text{mul}([e'_i])\right|\text{ for }i\in
\{1,\dotsc,n\}\\
e'_i\subseteq f_i\cap f^{-1}_{i-1}\text{ for }i\in \{1,\dotsc,n\}\\
e_i\text{ is the }i\text{th subedge of }e'_i\text{ relative to
}f_i\text{ for }i\in \{1,\dotsc,n\}
\end{gathered}
  \end{equation*}
We see that $f_n=f_0$, that $e'_n$ is the original edge of $f_0$
immediately following $e_0$, that $e_n$ is the terminal subedge of
$e'_n$ relative to $f_0$, that $e_0$ and $e_n$ are equivalent in an
orientation-preserving way, that $e_1$ is the terminal subedge of a
positive original edge of $f_0^{-1}$ and that none of the edges
$e_2,\dotsc,e_{n-1}$ is the terminal subedge of an original edge
relative to either face containing it.  Corresponding statements hold
if $e_0$ is an initial subedge of a negative original edge of $f_0$.

The previous paragraph implies for every face $f$ of $Q$ that the
terminal subedges of positive original edges of $f$ and $f^{-1}$ and
the initial subedges of negative original edges of $f$ and $f^{-1}$
are all equivalent and they are not equivalent to any other such edges
of other faces.  This and the results of the next-to-last paragraph
establish a bijection between the 1-cells of $M$ and pairs of faces of
$Q$.  Similarly, the last paragraph implies for every face of $Q$ that
its original vertices are equivalent.  This and the results of the
next-to-last paragraph imply that $M$ has just one 0-cell.

This proves Theorem~\ref{thm:main1}.
 
\end{proof}

We denote by $Q^*$ the subdivision of $P$ obtained by replacing the
multiplier function $\mul$ by $-\mul$.

\begin{thm}\label{thm:main2} Let $P$ be a faceted 3-ball, let $\ze$ be
an orientation-reversing face-pairing on $P$, and let $\mul$ be a
multiplier function for $(P,\ze)$.  Then the dual of the link of the
vertex of $M$ is isomorphic to $\partial Q^*$ in an
orientation-reversing way.
\end{thm}
\begin{proof} The proof is an adaptation of the arguments for the
analogous results in \cite{twist} and \cite{heegaard}.  Suppose $f$
is a face of $P$ and $e$ is an edge of $P$ in $f$. First suppose that
$e$ is a positive edge. Let $a$ be the initial vertex of $e$ relative
to $f$, let $b$ be the terminal vertex of $e$ relative to $f$, and let
$h$ be the edge of $Q$ preceding $e$ in $f$. Let $x$ be the vertex of
$M$. The image of $\link(a,Q)$ in $\link(x,M)$ has a vertex
corresponding to $h$, and this vertex is in a chain of $\ell([e])
\mul([e]) + 1$ faces; the first face is the image of $\link(a,Q)$, the
last face is the image of $\link(b,Q)$ and all of the other faces are
digons which are the images of links of vertices of $Q$ that are not
vertices of $P$. Similarly, if $e$ is a negative edge, $a$ is the
terminal vertex of $e$ relative to $f$, $b$ is the initial vertex of
$e$ relative to $f$, and $h$ is the edge of $Q$ following $e$ in $f$,
then the vertex corresponding to $h$ in the image of $\link(a,Q)$ in
$\link(x,M)$ is in a chain of $\ell([e])\left|\mul([e])\right|+1$
faces joining the images of $\link(a,Q)$ and $\link(b,Q)$. So in each
case, in the dual of $\link(x,M)$ there is a segment subdivided into
$\ell([e])\left|\mul([e])\right|$ edges which joins the duals of the
images of $\link(a,Q)$ and $\link(b,Q)$.

We next need to see how these segments fit together. We suppose for
convenience that $e$ is a positive edge. Let $e'$ be the edge of $P$
that precedes $e$ in $f$ and let $e''$ be the edge of $P$ that follows
$e$ in $f$.  If $e'$ is also a positive edge, then in the dual of
$\link(x,M)$ there is a face containing a pair of adjacent segments,
subdivided into $\ell([e]) \mul([e])$ and $\ell([e']) \mul([e'])$
edges. A similar statement holds if $e''$ is a positive edge. If $e'$
is a negative edge, then the edge of $Q$ preceding $e$ in $f$ is a
sticker, and is the same as the edge of $Q$ following $e'$ in
$f$. This sticker is the edge $h$ of the previous paragraph for both
$e$ and $e'$.  So in the dual of $\link(x,M)$ the segments
corresponding to $e$ and $e'$ are adjacent in some face.  If $e''$ is
a negative edge, then the terminal subedge of $e$ in $f$ and the
initial subedge of $e''$ in $f$ are equivalent to a sticker in the
face $f^{-1}$, and so there is a sticker in the dual of $\link(x,M)$
between the segments corresponding to $e$ and $e''$.  A similar
analysis holds if $e$ is a negative edge.

This implies that in the dual of $\link(x,M)$ there is a face
corresponding to $f$ that is cellularly homeomorphic to the face
corresponding to $f$ in $Q^*$.  This correspondence between faces of
$Q^*$ and faces of the dual of $\text{link}(x,M)$ respects adjacency
of faces.  So the dual of $\text{link}(x,M)$ is cellularly
homeomorphic to $\partial Q^*$. It follows as in \cite{twist} that
this homeomorphism reverses orientation.

\end{proof}

\begin{remark}\label{remark:main3} The proof of
Theorem~\ref{thm:main2} interpreted in terms of dual cap subdivision
shows just as in \cite{twist} and \cite{heegaard} that the manifolds
$M(P,\ze,\mul)$ and $M(P,\ze,-\mul)$ are homeomorphic by means of a
map which establishes a duality between these cell complexes.
\end{remark}

\section{Heegaard diagrams for bitwist manifolds}\label{sec:heegaard}\nosubsections

Let $M = M(P,\ze,\mul)$ be a bitwist manifold, let $Q$ be the
corresponding subdivision of $P$, and let $\zd$ be the corresponding
bitwisted face-pairing on $Q$. As in \cite[Section 4]{heegaard}, one
can construct the \textbf{edge pairing surface} $S$ of $(Q,\zd)$. For
each face $f$ in $Q$, there is a CW structure on a closed disk $F_f$
such that i) $F_f$ has a single 2-cell whose interior is the interior
of $F_f$, ii) there is a continuous cellular map $\varphi_f\co F_f\to
f$ whose restriction to each open cell in a homeomorphism, and iii)
there is a continuous cellular map $\psi_f \co F_f\to f^{-1}$ whose
restriction to each open cell is a homeomorphism. (And also
$\varphi_f$ and $\psi_f$ are compatible with respect to the
face-pairing.) Let $Y$ be the quotient of the union of the 1-skeleton
$X$ of $Q$ and the finite union of the complexes $\partial F_f \times
[0,1]$, one for each pair $(f,f^{-1})$, under the equivalence relation
generated by the identifications of $(x,0)$ with $\varphi_f(x)$ and
$(x,1)$ with $\psi_f(x)$ for $x\in \partial F_f$. Then $Y$ is an
orientable closed surface, and the dual cap subdivision of $Y$ is the
edge pairing surface $S$. (See \cite[Section 3]{heegaard} for the
definition of the dual cap subdivision. The dual cap subdivision of a
2-complex is obtained from its barycentric subdivision by removing the
edges joining vertices to barycenters of faces.) Edges of $S$ that are
contained in $X$ or disjoint from $X$ are called \textbf{vertical},
and the other edges of $S$ are called \textbf{diagonal}. Edges of $S$
that are not contained in edges of $Y$ are called \textbf{meridian}
edges, and edges of $S$ contained in edges of $Y$ are called
\textbf{nonmeridian} edges.

\begin{thm}\label{thm:heegaard}
Let $M = M(P,\ze,\mul)$ be a bitwist manifold, and let $S$ be
the edge pairing surface for the associated bitwisted face pairing.
The union $\cV$ of the vertical meridian edges is a basis of
meridian curves for $S$, and the union $\cD$ of the diagonal
meridian edges is a basis of meridian curves for $S$. Furthermore
$(S,\cV,\cD)$ is a Heegaard diagram for $M$.
\end{thm}
\begin{proof}
Since $M = Q/\zd$ is a manifold with a single vertex, this follows
immediately from \cite[Theorem 4.2.1]{heegaard}.

\end{proof}

\begin{figure}
\centerline{\includegraphics{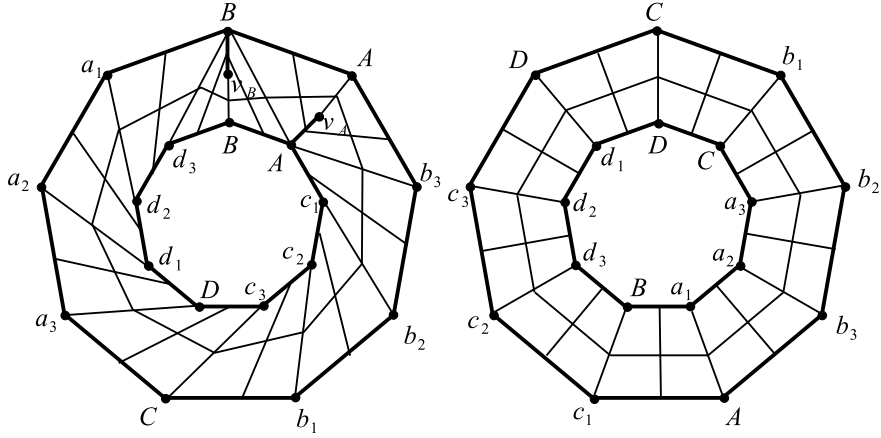}} \caption{$\partial
F_{f_1}\times [0,1]$ and $\partial F_{f_2}\times [0,1]$ for the
example from Section~\ref{sec:prelimex}.} \label{fig:surfaceex}
\end{figure}

Figure~\ref{fig:surfaceex} shows the union of $\partial
F_{f_1}\times [0,1]$ and $\partial F_{f_2}\times [0,1]$ for the
example from Section~\ref{sec:prelimex}, where $f_1$ is the
triangle $ABC$, $f_2$ is the triangle $ACD$, and the two sides of
the stickers have been identified.

\begin{figure}
\centerline{\includegraphics{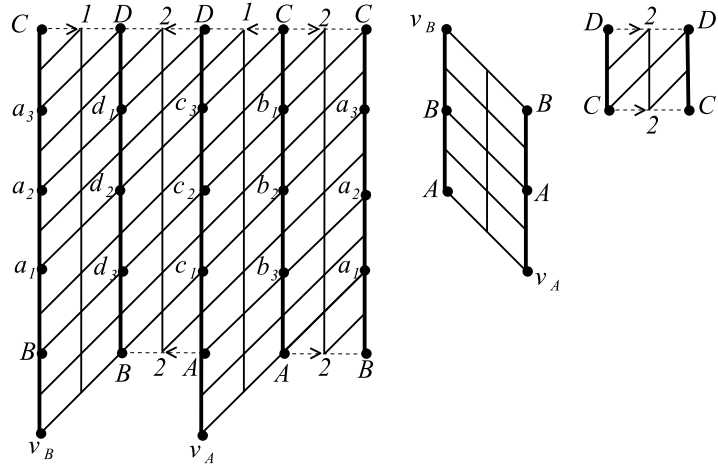}} \caption{The edge cycle
cylinders for the example from Section~\ref{sec:prelimex}.}
\label{fig:edgecylex}
\end{figure}

\begin{figure}
\centerline{\includegraphics{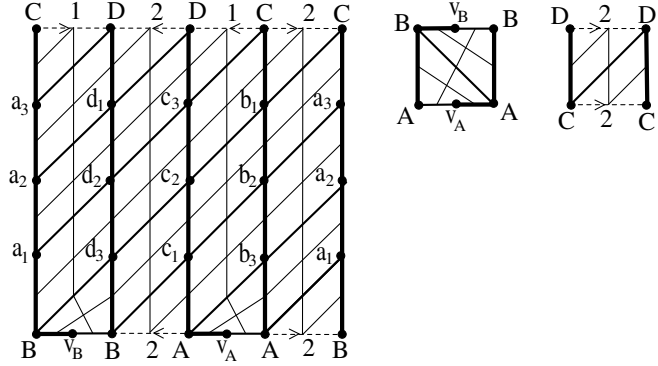}} \caption{Another view of
the edge cycle cylinders for the example from Section~\ref{sec:prelimex}.}
\label{fig:edgecyl2}
\end{figure}

As in \cite{heegaard}, the surface $S$ can also be decomposed into
edge cycle cylinders. The only difference from the construction in
\cite{heegaard} is that if $f$ is a face of $P$ and $e$ is either a
positive original edge which is preceded by a sticker or a negative
original edge that is followed by a sticker, then the sticker is
included with that edge in the construction of the edge cycle
cylinder. For example, Figure~\ref{fig:edgecylex} shows, for the
example from Section~\ref{sec:prelimex}, the edge cycle
cylinders. Figure~\ref{fig:edgecyl2} shows, for the same example, the
edge cycle cylinders with the stickers pushed back to be horizontal
edges. Note that, in this view, vertical meridian edges are drawn
vertically and diagonal meridian edges are drawn diagonally. This view
makes the effect of adding the stickers more apparent. When a diagonal
meridian edge crosses a sticker, it changes direction. This reflects
the difference in directions of twists corresponding to edge
cycles with positive multipliers and edge cycles with negative
multipliers.

Let $C$ be an edge cycle cylinder, where as in
Figure~\ref{fig:edgecyl2} we have pushed the stickers back to be
horizontal. Let $\alpha$ (resp. $\alpha'$)  be a minimal union of
vertical (resp. diagonal) meridian edges that joins the two
horizontal ends of $C$, chosen so that $\partial \alpha = \partial
\alpha'$. Let $\beta$ be a simple closed curve in $C$ that separates
the ends of $C$, and let $m = \mul(E)$, where $E$ is the edge
cycle associated to $C$.  Then $\alpha'$ is isotopic rel endpoints
to $\tau^m(\alpha)$, where $\tau$ is a Dehn twist along $\beta$.
Furthermore, as one repeats this construction for the other edge
cycle cylinders, the directions of the Dehn twists can all be chosen
consistently with respect to an orientation of $S$.

\begin{thm}\label{thm:dehntwist} Let $M = M(P,\ze,\mul)$ be a bitwist
manifold, let $S$ be the edge pairing surface for the associated
bitwisted face pairing, and let $\cV = \{\alpha_1,\dots,\alpha_n\}$ be
the vertical meridian curves as in Theorem~\ref{thm:heegaard}. Let
$E_1,\dots,E_m$ be the edge cycles of $\ze$. For each
$i\in\{1,\dots,m\}$ let $C_i$ be the edge cycle cylinder associated to
$E_i$ and let $\tau_i$ be a Dehn twist along a simple closed curve in
$C_i$ which separates the ends of $C_i$. We choose the $\tau_i$'s so
that they twist in consistent directions with respect to a fixed
orientation of $S$. Let $\tau = \tau_1^{\mul(E_1)} \circ \dots \circ
\tau_m^{\mul(E_m)}$. Then
$(S,\cV,\{\tau(\alpha_1),\dots,\tau(\alpha_n)\})$ is a Heegaard
diagram for $M$.
\end{thm}
\begin{proof}
This follows immediately from Theorem~\ref{thm:heegaard} and the
discussion in the paragraph before the statement of the theorem.

\end{proof}

The construction of corridor complex links for bitwist 3-manifolds is
the same as their construction in \cite[Section 6]{heegaard} for
twisted face-pairing manifolds, though the framings change because of
the signs of the multipliers. We first recall the construction of
corridor complex links.

Suppose $P$ is a faceted 3-ball, $\ze$ is an orientation-reversing
face-pairing on $P$, and $\mul$ is a multiplier function for $\ze$.
Let $M = M(P,\ze,\mul)$ be the associated bitwist 3-manifold. We form
a \textbf{corridor complex} for $\ze$ as follows. We choose a pair $f_1$
and $f_2$ of faces in $\partial P$ that are matched by $\ze$, and
choose an edge-path arc in the 1-skeleton of $\partial P$ that joins a
corner of $f_1$ to its image under $\ze$ in $f_2$. We then split this
edge-path to a thin corridor. This gives a new cell structure on
$\partial P$ in which the old faces $f_1$ and $f_2$ have been joined
by the corridor into a single face. We do this successively for all of
the face pairs of $\partial P$, and call the resulting cell structure
on $\partial P$ the corridor complex $C$.

\begin{figure}
\centerline{\includegraphics{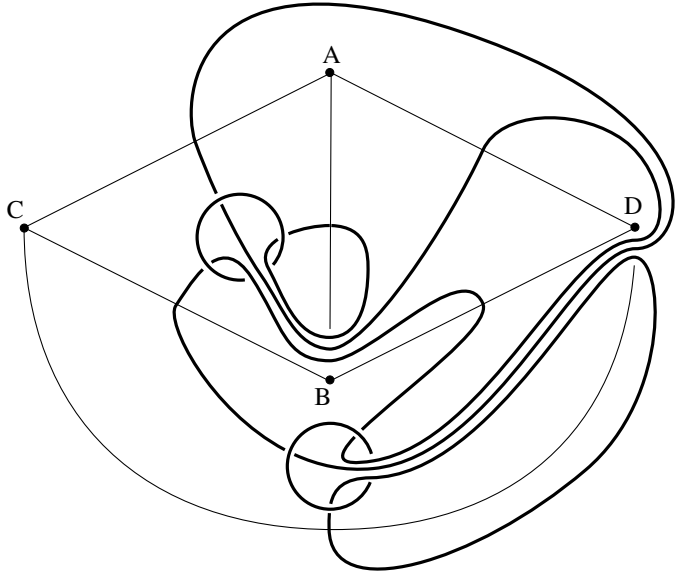}} \caption{A corridor
complex link for the example from Section~\ref{sec:prelimex}.}
\label{fig:corlnkex}
\end{figure}

We next describe a link $L$ in $S^3$ in terms of its projection to
$C$. For each face of $C$ there is an unknotted component of $L$ that
lies in one of the old faces that are part of that face; we call this
component a \textbf{face component}. Next consider one of the old faces
$f$ that contains a face component. Each edge of that old face
corresponds to an edge of the corresponding face in the corridor
complex. For each such edge $e$, $L$ contains an arc which enters the
old face from the barycenter of the edge, crosses under the face
component in the old face, crosses over the face component, goes
through the corridor, and ends at the barycenter of the edge
$\ze_f(e)$. These arcs are constructed so that they have no
self-crossings or intersections with other such arcs from that
face. We construct these arcs for each face of the corridor complex.
Suppose $e$ is one of the original edges in $P$. If $e$ has not been
split in the construction of the corridor complex, then at the
barycenter of $e$ we have the ends of the arcs from the two faces that
contain $e$ (or from the face that meets $e$ with multiplicity
two). If $e$ has been split in the construction of the corridor
complex, then we join the ends of the two corresponding arcs by an
arc that goes under the arcs in the corridor. The union of all of
these arcs is a finite set of components of $L$ that are called \textbf{
edge components}. Each edge component crosses exactly those edges of
$C$ which correspond to an edge cycle of $\ze$. The \textbf{corridor
complex link} $L$ is the union of the face components and the edge
components. We call $L$ a corridor complex link for $(P,\ze)$. A
corridor complex link for the example from Section~\ref{sec:prelimex}
is shown in Figure~\ref{fig:corlnkex}.

\begin{thm}\label{thm:bilink}
Let $M = M(P,\ze,\mul)$ be a bitwist 3-manifold, and let $L$ be
the corresponding corridor complex link. Define a framing on $L$ by
giving each face component framing $0$ and giving the edge component
corresponding to an edge cycle $E$ the framing $\mul(E)^{-1}$
plus the blackboard framing of the edge component. Then Dehn surgery
on the framed link $L$ yields $M$.
\end{thm}
\begin{proof}This follows easily from the proofs of
\cite[Theorem 6.2.2]{heegaard} and \cite[Theorem 6.1.2]{heegaard}.
The proof of \cite[Theorem 6.2.2]{heegaard} goes through in this
greater generality until the last paragraph, when it refers to
\cite[Theorem 6.1.2]{heegaard}. The statement and proof of
\cite[Theorem 6.1.2]{heegaard} go through in this greater generality.

\end{proof}

\section{Generalizing framings of corridor complex links}
\label{sec:framings}\subsections

In this section we develop some of the machinery needed for the
proof of Theorem~\ref{thm:allm}. We first discuss some well-known
techniques for changing framed surgery descriptions of
3-manifolds. We then show that, in a sense made precise in
Theorem~\ref{thm:linksum}, connected sums of corridor complex
links are corridor complex links. Theorem~\ref{thm:consum}, that
connected sums of bitwist manifolds are bitwist manifolds, follows
easily. We next consider a special family of face-pairings called
reflection face-pairings, and use them to show that every lens
space is a twisted face-pairing 3-manifold. This allows us to
prove Theorem~\ref{thm:glinks}, which states that if $L$ is a
complex corridor link, then for any choices of framings for the
edge components we still get a bitwist manifold by framed surgery.

\subsection{Dehn surgery preliminaries}\label{sec:prelim}
\nosubsubsections

We collect some well-known facts about Dehn surgery which will be
used later.

We first discuss Rolfsen twists.  They appear on page 162 of
\cite{GS}, they appear in Sections 16.4, 16.5 and 19.4 of \cite{PS} as
Fenn-Rourke moves, and they appear in Section 9.H of \cite{R}.  For
this let $L$ be a link in $S^3$ framed by the elements of $\bQ\cup
\{\infty\}$.  Let $J$ be an unknotted component of $L$.  Then
$L\setminus J$ is contained in a closed solid torus $T$, which is the
complement in $S^3$ of a regular neighborhood of $J$.  Let $\zt$ be a
right hand Dehn twist of $T$.  Let $n\in \bZ $.  Let $L'$ be the link
gotten from $L$ by applying $\zt^n$ to $L\setminus J$.  We frame $L'$
as follows.  If the $L$-framing of $J$ is $r$, then the $L'$-framing
of $J$ is $\frac{1}{n+\frac{1}{r}}$.  If $K$ is a component of $L$
other than $J$ with framing $r$, then the image of $K$ in $L'$ has
framing $r+n\cdot\textrm{lk}^2(J,K)$, where $\textrm{lk}(J,K)$ is the
linking number of $J$ and $K$ after orienting $J$ and $K$ arbitrarily.
When $n=1$, we say that $L'$ is obtained from $L$ by performing a
Rolfsen twist about $J$.  In general we obtain $L'$ by performing $n$
Rolfsen twists about $J$.  We are interested in Rolfsen twists because
the manifold obtained by Dehn surgery on $L'$ is homeomorphic to the
manifold obtained by Dehn surgery on $L$.

We next discuss slam-dunks.  These appear on page 163 of \cite{GS}.
Let $L$ be a framed link in $S^3$.  Suppose that one component $K$
of $L$ is a meridian of another component $J$ and that $K$ is
contained in a topological ball in $S^3$ which meets no components of
$L$ other than $J$ and $K$.  Suppose that the framing of $J$ is
$n\in \bZ $ and that the framing of $K$ is $r\in \bQ\cup\{\infty\}$.
Let $L'$ be the framed link obtained from $L$ by deleting $K$ and
changing the framing of $J$ to $n-\frac{1}{r}$.  We say that $L'$ is
obtained from $L$ by performing the slam-dunk which removes $K$.
The manifold obtained by Dehn surgery on $L'$ is homeomorphic to the
manifold obtained by Dehn surgery on $L$.

\subsection{Connected sums of corridor complex links}
\label{sec:linksum}\nosubsubsections

Here we establish the fact that the links
obtained from the corridor construction are
closed under the operation of connected sum in a certain restricted
sense.

\begin{figure}
\centerline{\includegraphics{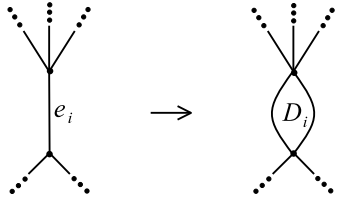}} \caption{Replacing $e_i$
with a digon $D_i$.} \label{fig:replace}
\end{figure}

We begin with two faceted 3-balls $P_1$ and $P_2$.  For $i\in\{1,2\}$
let $\ze_i$ be an orientation-reversing face-pairing on $P_i$ with
multiplier function $\mul_i$, and let $M_i=M(P_i,\ze_i,\mul_i)$. For
$i\in\{1,2\}$ let $L_i$ be the link corresponding to $M_i$ as in
Theorem~\ref{thm:bilink}.  For $i\in\{1,2\}$, let $C_i$ be an edge
component of $L_i$ and let $e_i$ be an edge of $P_i$ which lies in the
$\ze_i$-edge cycle corresponding to $C_i$. We assume that either $e_1$
has distinct vertices or $e_2$ has distinct vertices. Let $P'_i$ be
the faceted 3-ball obtained from $P_i$ by replacing $e_i$ with a digon
$D_i$ for $i\in\{1,2\}$.  See Figure~\ref{fig:replace}.  Because
either $e_1$ has distinct vertices or $e_2$ has distinct vertices, we
obtain a faceted 3-ball $P$ from $P'_1$ and $P'_2$ by cellularly
identifying $D_1$ and $D_2$.  We refer to $P$ as a \textbf{connected
sum} of $P_1$ and $P_2$ along $e_1$ and $e_2$. The face-pairings
$\ze_1$ and $\ze_2$ induce a face-pairing $\ze$ on $P$. Except for
choices to be made involving corridors along either $e_1$ or $e_2$,
the corridor constructions for $(P_1,\ze_1)$ and $(P_2,\ze_2)$ which
give rise to $L_1$ and $L_2$ induce a corridor construction for
$(P,\ze)$, which gives rise to an unframed link $L$.  The isotopy type
of $L$ is uniquely determined by $L_1$, $L_2$ and the identification
of $D_1$ and $D_2$. It is easy to see that $L$ is a connected sum of
$L_1$ and $L_2$ which joins $C_1$ and $C_2$.  We summarize this
paragraph in the following theorem.

\begin{thm}[]\label{thm:linksum} Let $P_1$ and $P_2$ be faceted
3-balls with orientation-reversing face-pairings $\ze_1$ and $\ze_2$.
Let $L_1$ and $L_2$ be corresponding unframed corridor complex links.
Let $C_1$ be an edge component of $L_1$, and let $C_2$ be an edge
component of $L_2$.  Let $e_1$ be an edge of $P_1$ which lies in the
$\ze_1$-edge cycle corresponding to $C_1$, and let $e_2$ be an edge of
$P_2$ which lies in the $\ze_2$-edge cycle corresponding to
$C_2$. Suppose that either $e_1$ has distinct vertices or $e_2$ has
distinct vertices.  Let $P$ be a connected sum of $P_1$ and $P_2$
along $e_1$ and $e_2$, and let $L$ be the corresponding connected sum
of $L_1$ and $L_2$ which joins $C_1$ and $C_2$.  Then $L$ is an
unframed corridor complex link associated to the orientation-reversing
face-pairing on $P$ induced by $\ze_1$ and $\ze_2$.
\end{thm}
\begin{proof} This is clear from the previous paragraph.

\end{proof}

Suppose $P_1$ and $P_2$ are faceted 3-balls. For $i\in\{1,2\}$ let
$\ze_i$ be an orientation-reversing face-pairing on $P_i$ and let
$\mul_i$ be a multiplier function for $\ze_i$.  Let $e_1$ be an edge
in $P_1$ and let $e_2$ be an edge in $P_2$ such that $\mul_1([e_1]) =
\mul_2([e_2])$.  Then the multiplier functions $\mul_1$ and $\mul_2$
induce a multiplier function for the face-pairing induced by $\ze_1$
and $\ze_2$ on the connected sum of $P_1$ and $P_2$ along $e_1$ and
$e_2$.

\subsection{Connected sums of bitwist manifolds}
\label{sec:consum}\nosubsubsections

\begin{thm}[]\label{thm:consum} The connected sum of two
bitwist manifolds is a bitwist manifold.
\end{thm}
\begin{proof} Let $P$ be the faceted 3-ball with just two faces which
are degenerate pentagons as in Figure~\ref{fig:consump}.  Let $\ze$ be
the face-pairing on $P$ which fixes the edge common to the two faces,
and let mul be the multiplier function for $\ze$ indicated in
Figure~\ref{fig:consump}.  Figure~\ref{fig:consuml} shows a corridor
complex and a corridor complex framed link $L$ for $\ze$ and mul.

\begin{figure}
\centerline{\includegraphics{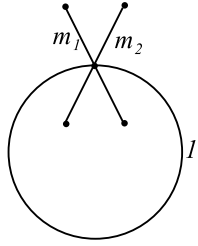}} \caption{The faceted
3-ball $P$ and edge cycle multipliers.}
\label{fig:consump}
\end{figure}

\begin{figure}
\centerline{\includegraphics{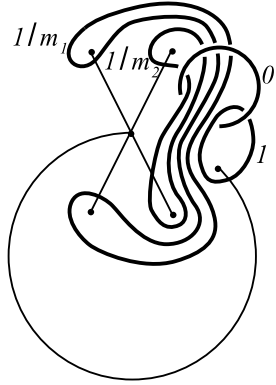}} \caption{The framed
corridor complex link $L$.}
\label{fig:consuml}
\end{figure}

Now let $P_1$ and $P_2$ be faceted 3-balls with face-pairings and
multiplier functions which give rise to bitwist manifolds $M_1$ and
$M_2$.  We choose one of the two edges of $P$ in the $\ze$-edge cycle
with multiplier $m_1$, and we form a connected sum $P'_1$ of $P$ and
$P_1$ along this edge and any edge of $P_1$.  Next we choose one of
the two edges of $P$ in the $\ze$-edge cycle with multiplier $m_2$.
This edge corresponds to an edge of $P' _1$.  We form a connected sum
$P'_2$ of $P'_1$ and $P_2$ along this edge and any edge of $P_2$.
Theorem~\ref{thm:linksum} easily implies that we obtain a twisted
face-pairing manifold $M$ which is the connected sum of $M_1$, $M_2$,
and a manifold which is obtained by Dehn surgery on a framed link
which consists of two simply linked unknots with framings 0 and 1.
This third connected summand is the 3-sphere.  Thus $M$ is the
connected sum of $M_1$ and $M_2$.

This proves Theorem~\ref{thm:consum}.

\end{proof}

\subsection{Reflection face-pairings}\label{sec:reflection}
\nosubsubsections

We next consider face-pairings of a very special sort.
We assume that our model faceted 3-ball $P$ can be identified with
the closed unit ball in $\bR^3$ so that the following holds.  The
intersection of the unit sphere with the $xy$-plane is a union of
edges of $P$ and the model face-pairing $\ze$ on $P$ is given by
reflection in the $xy$-plane.  In other words, we have cell
structures on both the northern and southern hemispheres of the unit
sphere in $\bR^3$, and the face-pairing maps of the model
face-pairing $\ze$ are given by the map $(x,y,z)\mapsto (x,y,-z)$,
which is therefore a cellular automorphism of $P$.  In this case we
call $P$ a \textbf{reflection faceted 3-ball}, and we call $\ze$ a
\textbf{reflection face-pairing}.  Using the identification of $P$
with the closed unit ball in $\bR^3$, we speak of the
\textbf{equator} of $P$ and the \textbf{northern} and
\textbf{southern hemispheres} of $P$.

\begin{figure}
\centerline{\includegraphics{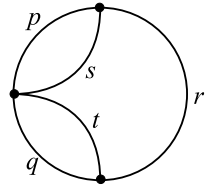}} \caption{The diagram
corresponding to $P$, $\ze$, and $\mul$.} \label{fig:linkexp}
\end{figure}

Let $P$ be a reflection faceted 3-ball with reflection face-pairing
$\ze$ and multiplier function $\mul$. As in
Figure~\ref{fig:linkexp}, we can describe $P$, $\ze$, and $\mul$
using a diagram which consists of a cellular decomposition of
a closed disk together with a nonzero integer for every edge.  We
view this closed disk as the northern hemisphere of $P$.  Hence we
have the cellular decomposition of the northern hemisphere of $P$,
which therefore determines the cellular decomposition of the
southern hemisphere of $P$, and the integer attached to the edge $e$
is $\mul([e])$.  We sometimes allow ourselves the liberty of
attaching 0 to an edge instead of a nonzero integer.  Attaching 0 to
an edge means that every edge in the corresponding $\ze$-edge cycle
collapses to a vertex.

Let $P$ be a reflection faceted 3-ball with reflection face-pairing
$\ze$.  Suppose given a multiplier function $\mul$
for $\ze$, and let $M$ be the associated bitwist manifold.
Theorem~\ref{thm:bilink} describes a framed link in the 3-sphere
$S^3$ such that Dehn surgery on this framed link gives $M$. In this
paragraph we describe another framed link $L$ in $S^3$ such that
Dehn surgery on $L$ also gives $M$.  We construct $L$ as follows. We
identify $P$ with the closed unit ball in $\bR^3$ as in the
definition of reflection faceted 3-ball.  For every edge $e$ of the
northern hemisphere of $P$ we choose an open topological ball
$B_e\subseteq \bR^3$ such that $B_e\cap
\partial P$ is a topological disk which meets $e$ and is disjoint from
every edge of $P$ other than $e$.  We assume that such topological
balls corresponding to distinct edges are disjoint.  For every
face $f$ of the northern hemisphere of $P$ we construct an unknot
$C_f$ in the interior of $f$ such that if $e$ is an edge of $f$,
then $C_f$ meets $B_e$. These unknots are all components of $L$
with framings 0.  We call these components of $L$ \textbf{face
components}.  Let $\zs\in\{\pm 1\}$. Every edge $e$ of $P$ in the
northern hemisphere also gives a component $C_e$ of $L$, called an
\textbf{edge component}, as follows.  Let $e$ be an edge in the
equator of $P$ contained in the face $f$ of the northern
hemisphere. The $\ze$-edge cycle of $e$ is just $\{e\}$.  We
define $C_e$ to be a meridian of $C_f$ contained in $B_e$ with
framing $\zs/\mul(\{e\})$.  Now let $e$ be
an edge of the northern hemisphere of $P$ not contained in the
equator.  Let $f$ and $g$ be the faces of $P$ which contain $e$.
Let $x$ be a point of $f\cap B_e$ separated by $C_f$ from
$\partial f$, and let $y\ne x$ be a point of $g\cap B_e$ separated
by $C_g$ from $\partial g$. The $\ze$-edge cycle of $e$ is $E =
\{e,\ze_f(e)\}$.  We define $C_e$ to be an unknot in $B_e$ with
framing $\zs/\mul(E)$ such that $P\cap C_e$ is a properly
embedded arc in $P\cap B_e$ joining $x$ and $y$. This defines $L$.

\begin{figure}
\centerline{\includegraphics{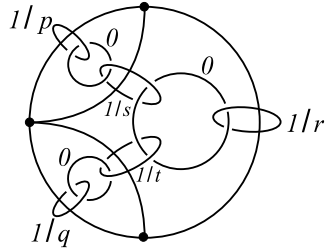}} \caption{The framed link
$L$.} \label{fig:linkexl}
\end{figure}

\begin{ex}[]\label{ex:linkex} Let $P$ be the reflection faceted 3-ball
with reflection face-pairing, and multiplier function given by the
diagram in Figure~\ref{fig:linkexp}.  Figure~\ref{fig:linkexl} shows
the framed link $L$ constructed above from these data using $\zs=1$.
\end{ex}

\begin{thm}[]\label{thm:link} Let $P$ be a reflection faceted 3-ball
with reflection face-pairing $\ze$.  Suppose given a multiplier
function for $\ze$, and let $M$ be the associated bitwist manifold.
Let $L$ be the framed link in $S^3$ constructed above.  Then Dehn
surgery on $L$ gives $M$.
\end{thm}
\begin{proof} Since $L$ is amphicheiral, multiplying all framings by
$-1$ does not change the resulting manifold.  So we may assume that
$\zs=1$.  We show how to adapt \cite[Theorem 6.1.2]{heegaard} to the
present situation.

We construct a handlebody $H$ as follows.  We still identify $P$ with
the closed unit ball in $\bR^3$.  Let $B$ be the topological ball
which is the closure in $S^3$ of $S^3\setminus P$.  We construct $H$
by attaching handles to $B$ as follows.  Let $f$ and $f^{-1}$ be faces
of $P$ paired by $\ze$.  Then $f$ and $f^{-1}$ are joined by a
vertical circular cylinder.  We attach such a cylinder to $B$.  Doing
this for every pair of faces of $P$ yields our handlebody $H$.  It is
clear that the closure in $S^3$ of $S^3\setminus H$ is also a
handlebody.  We identify the components of $L$ with curves in
$\partial H$ in a straightforward way.

As in \cite[Theorem 6.1.2]{heegaard}, let $S$ be the edge pairing
surface for the bitwisted face-pairing $\zd$, let $\za_1,\dotsc,\za_n$
be the vertical meridian curves of $S$ and let $\zb_1,\dotsc,\zb_m$ be
core curves for the edge cycle cylinders.  Then there exists a
homeomorphism $\zv\co S\to \partial H$ such that $\zv(\za_i)$ is the
face component of $L$ corresponding to $\za_i$, this face component
being a meridian of $H$, for every $i\in\{1,\ldots,n\}$.  We also have
that the edge components of $L$ are parallel copies of
$\zv(\zb_1),\dotsc,\zv(\zb_m)$.  The framing determined by $\partial
H$ of every edge component of $L$ is 0. Just as in the proof of
Theorem~\ref{thm:bilink}, the statement and proof of \cite[Theorem
6.1.2]{heegaard} go through in this greater generality.  So Dehn
surgery on $L$ gives $M$.

\end{proof}

\subsection{Lens spaces}\label{sec:lens}\nosubsubsections

In this subsection we show that every lens space is a twisted
face-pairing manifold. We will use this in the proof of
Theorem~\ref{thm:glinks}.

\begin{figure}
\centerline{\includegraphics{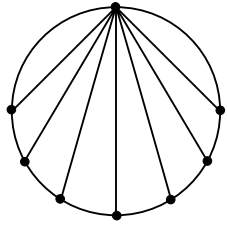}} \caption{Top view of a
scallop.} \label{fig:scallop}
\end{figure}

We begin by defining the notion of a scallop.  A \textbf{scallop} is a
reflection faceted 3-ball $P$ (defined in
Section~\ref{sec:reflection}) whose northern hemisphere has a cell
structure essentially as indicated in Figure~\ref{fig:scallop}.  More
precisely, every vertex of a scallop $P$ lies on the equator of $P$,
$P$ contains a vertex $v$ such that every edge of $P$ not contained in
the equator of $P$ joins $v$ with another vertex, and every vertex of
$P$ other than $v$ is joined with $v$ by at least one edge.  So the
northern hemisphere of a scallop might consist of just a monogon.
Otherwise it is subdivided into digons and triangles, in which case it
has at least two digons, but it may have arbitrarily many digons.

\begin{figure}
\centerline{\includegraphics{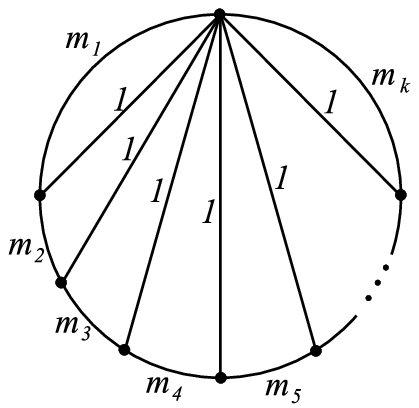}} \caption{The diagram for
$P$, $\ze$ and mul.} \label{fig:lensp}
\end{figure}

\begin{thm}[]\label{thm:lens} Let $P$ be a scallop with $k$ faces
in its northern hemisphere.  Let $\ze$ be a reflection face-pairing
on $P$, let mul be a multiplier function for $\ze$, and let
$M=M(P,\ze,\text{mul})$.  Suppose that $P$, $\ze$, and mul are given
by the diagram in Figure~\ref{fig:lensp}, where $m_1>0$, $m_k>0$,
and $m_i\ge 0$ for $i\in\{2,\dotsc,k-1\}$.  (If a multiplier is 0,
then the corresponding edge in Figure~\ref{fig:lensp} collapses to a
vertex of $P$.)  Define integers $a_1,\dotsc,a_k$ so that $a_1=m_1$
if $k=1$ and if $k>1$, then $a_1=m_1+1$, $a_k=m_k+1$, and
$a_i=m_i+2$ for $i\in\{2,\dotsc,k-1\}$.  Then there exist relatively
prime positive integers $p\ge q$ such that $M$ is homeomorphic to
the lens space $L(p,q)$, where
  \begin{equation*}
\frac{p}{q} = [a_1,-a_2,a_3,\dots,(-1)^{k+1}a_k] =
a_1-\cfrac{1}{a_2-\cfrac{1}{a_3-\cdots-\cfrac{1}{a_{k-1}-\cfrac{1}{a_k}}}}.
\end{equation*}
(It is possible that $p=q=1$, in which case we obtain the 3-sphere.)
Furthermore, given relatively prime positive integers $p$ and $q$
with $p\ge q$, then there exists a unique sequence of integers
$m_1,\dotsc,m_k$ as above such that the above continued fraction
equals $p/q$.
\end{thm}

\begin{proof} Theorem~\ref{thm:link} implies that $M$ is given by Dehn
surgery on the framed link in Figure~\ref{fig:lensl}, where for
convenience we have chosen $\sigma = -1$. We repeat that if $m_i=0$
for some $i\in\{2,\dotsc,k-1\}$, then the corresponding edge in
Figure~\ref{fig:lensp} collapses to a vertex of $P$.  In this case
the corresponding component of the link in Figure~\ref{fig:lensl} is
to be removed.  This is consistent with the fact that any component
with framing $\infty$ may be removed from a framed link without
changing the resulting manifold.  We next use Kirby calculus to
simplify the framed link in Figure~\ref{fig:lensl}.  For every
$i\in\{1,\dotsc,k\}$ we perform the slam-dunk which removes the
component with framing $-1/m_i$. In doing this, the component linked
with the given component acquires the framing $m_i$.  We next
perform a Rolfsen twist about every component shown in
Figure~\ref{fig:lensl} with framing $-1$.  Every such component is
then removed, and 1 is added to the framing of the components linked
with it.  The resulting framed link is shown in Figure~\ref{fig:chain}.
It follows from page 272 of \cite{R} or page 108 of \cite{PS} or just
by iterating slam-dunks that $M$ is the lens space as stated in
Theorem~\ref{thm:lens}.

\begin{figure}
\centerline{\includegraphics{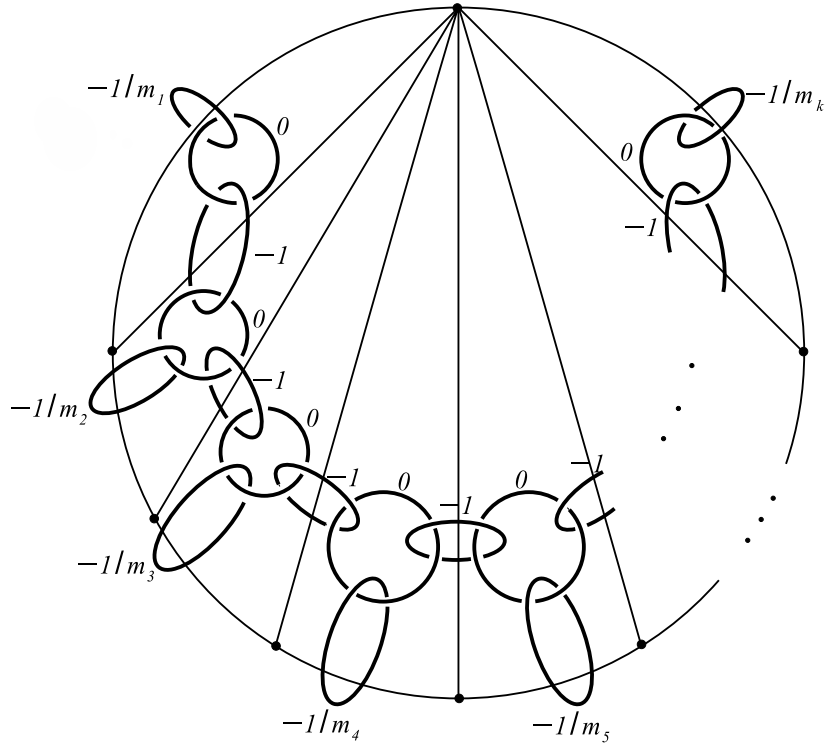}} \caption{The framed link
corresponding to Figure~\ref{fig:lensp}.} \label{fig:lensl}
\end{figure}

\begin{figure}
\centerline{\includegraphics{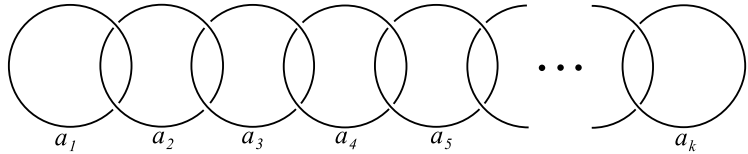}} \caption{Dehn surgery on
this framed link gives $M$.} \label{fig:chain}
\end{figure}

The uniqueness statement is well known.  For this, first note that
if $k=1$, then $a_1$ is an arbitrary positive integer.  If $k>1$,
then $a_1,\dotsc,a_k$ are arbitrary integers with $a_i\ge 2$ for
$i\in\{1,\dotsc,k\}$.  Given $p$ and $q$, we calculate
$a_1,\dotsc,a_k$ by modifying the division algorithm usually used to
calculate continued fractions.  Instead of taking the greatest
integer less than or equal to our given number, we take the least
integer greater than or equal to our given number.  The details are
left to the reader.

This proves Theorem~\ref{thm:lens}.

\end{proof}

\begin{cor}[]\label{cor:lens} Every lens space is a twisted
face-pairing manifold.
\end{cor}

\subsection{Changing the framings}
\label{sec:glinks}\nosubsubsections

Suppose given an orientation-reversing face pairing $\ze$ on a faceted
3-ball $P$. In Section~\ref{sec:heegaard} we construct a corridor
complex link $L$ by means of link projections. The face components of
$L$ correspond to the face-pairs of $\ze$, and the edge components of
$L$ correspond to the edge cycles of $\ze$. Given the extra
information of a multiplier function $\mul$, we define framings on the
components of $L$. We define the framing of each face component to be
0. If $C$ is an edge component, then we define the framing of $C$ to
be the blackboard framing of $C$ plus $\mul(E)^{-1}$, where $E$ is the
edge cycle corresponding to $C$.  By Theorem~\ref{thm:bilink},
performing Dehn surgery on $L$ with this framing obtains our bitwist
manifold $M(P,\ze,\mul)$. The following theorem states that if we
redefine the framing of $L$ by replacing each framing of an edge
component by an arbitrary rational number, then Dehn surgery on $L$
still obtains a bitwist manifold (usually constructed from a different
faceted 3-ball).

\begin{thm}[]\label{thm:glinks} Let $L$ be an unframed corridor
complex link.  We frame $L$ as follows.  Let $C$ be a component of
$L$. If $C$ is a face component, then we define the framing of $C$
to be 0. If $C$ is an edge component, then we define the framing of
$C$ to be an arbitrary rational number. Then Dehn surgery on $L$
with this framing obtains a bitwist manifold.
\end{thm}
\begin{proof} Let $P$ be a faceted 3-ball and let $\ze$ be an
orientation-reversing face pairing on $P$ such that $L$ is a
corridor complex link for $(P,\ze)$. Let $E_1,\dots,E_m$ be the edge
cycles, and let $C_1,\dots,C_m$ be the corresponding edge components
of $L$. For $i\in\{1,\dots,m\}$, let $b_i$ be the blackboard framing
of $C_i$ and let $\alpha_i \in \bQ$ such that $b_i + \alpha_i$ is
the framing on $C_i$. Let $\cN = \{i\in \{1,\dots,m\}\co \alpha_i\
\mbox{is not the reciprocal of an integer}\}$.

Suppose given $i\in\{1,\dots,m\}$. If $i\notin \cN$, then we define
the multiplier of $E_i$ to be $\mul(E_i) = 1/\alpha_i$.  If $i\in
\cN$, we in effect change the framing of $C_i$ by ``attaching a
scallop'' to our model faceted 3-ball, proceeding as follows.

Suppose that $i\in \cN$ and $\alpha_i \le 0$. Let $r_i = 1/(1 -
\alpha_i)$, and let $a_1,\dots,a_k$ be positive integers with $a_j \ge
2$ if $k > 1$ such that $r_i = [a_1,-a_2,a_3,\dots,(-1)^{k+1}a_k]$.
As in Theorem~\ref{thm:lens}, define $m_1,\dots,m_k$ by $m_1 = a_1$ if
$k=1$ and, if $k > 1$, $m_1 = a_1-1$, $m_k = a_k - 1$, and $m_j = a_j
- 2$ for $j\in \{2,\dots,k-1\}$. Let $P_i$ be the reflection faceted
3-ball shown in Figure~\ref{fig:lensp2}, and let $\ze_i$ be the
associated reflection face-paring. Define the multiplier of $E_i$ to
be $\mul(E_i) = 1$, and define the multiplier function on $(P_i,\ze_i)$ as
indicated in Figure~\ref{fig:lensp2}.

\begin{figure}
\centerline{\includegraphics{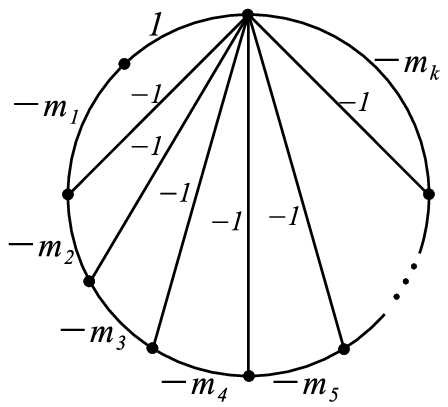}} \caption{The reflection
faceted 3-ball $P_i$ when $i\in \cN$ and $\alpha_i \le 0$.}
\label{fig:lensp2}
\end{figure}

Now suppose that $i\in \cN$ and $\alpha_i > 0$. Let $r_i = 1/(1 +
\alpha_i)$, and let $a_1,\dots,a_k$ be positive integers with $a_j
\ge 2$ if $k > 1$ such that $r_i =[a_1,-a_2,\dots,(-1)^{k+1}a_k]$.
As in Theorem~\ref{thm:lens}, define $m_1,\dots,m_k$ by $m_1 =
a_1$ if $k=1$ and, if $k > 1$, $m_1 = a_1-1$, $m_k = a_k - 1$,
and $m_j = a_j - 2$ for $j\in \{2,\dots,k-1\}$. Let $P_i$ be the
reflection faceted 3-ball shown in Figure~\ref{fig:lensp3}, and let
$\ze_i$ be the associated reflection face-paring. Define the
multiplier of $E_i$ to be $\mul(E_i) = -1$, and define the multiplier
function on $(P_i,\ze_i)$ as indicated in Figure~\ref{fig:lensp3}.

\begin{figure}
\centerline{\includegraphics{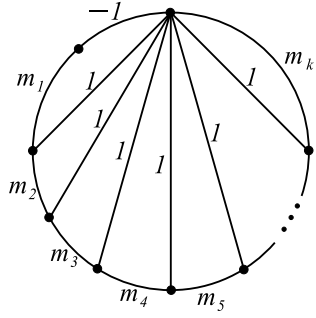}} \caption{The reflection
faceted 3-ball $P_i$ when $i\in \cN$ and $\alpha_i > 0$.} \label{fig:lensp3}
\end{figure}

We now construct the faceted 3-ball $P'$ and orientation-reversing
face-pairing $\ze'$ by repeated connect sums of $P$ with the faceted
3-balls $P_i$ for which $i\in \cN$. For each $i\in\cN$, we do this via
an edge in the edge cycle corresponding to $C_i$ and the edge in $P_i$
which is immediately to the left of the top vertex in
Figure~\ref{fig:lensp2} or \ref{fig:lensp3}. Since the multipliers are
compatible on edge cycles that are amalgamated, they define a
multiplier function for $\ze'$.

We next construct a framed corridor complex link for $(P',\ze')$. If
$i\in \cN$ and $\alpha_i \le 0$, then the link $K_i$ shown in
Figure~\ref{fig:glinks0} is a framed link for $(P_i,\ze_i)$ as in
Figure~\ref{fig:lensl}.  This framed link is in fact isotopic to a
framed corridor complex link for $(P_i,\ze_i)$.  If $i\in\cN$ and
$\alpha_i > 0$, then one gets a framed corridor complex link $K_i$
from the link in Figure~\ref{fig:glinks0} by multiplying the framing
of each component by $-1$. By repeated applications of
Theorem~\ref{thm:linksum}, one gets a framed corridor complex link $J$
for $(P',\ze')$.

Suppose $i\in\cN$ and $\alpha_i \le 0$. Figure~\ref{fig:linksum}
shows part of $J$ corresponding to $K_i$. As in the proof of
Theorem~\ref{thm:lens}, we can simplify this to obtain the framed
link in Figure~\ref{fig:glinks1}.  Again as in the proof of
Theorem~\ref{thm:lens}, by performing $k-1$ slam-dunks, we may
reduce $J$ to the framed link in Figure~\ref{fig:glinks2}. A similar
argument holds if $i\in\cN$ and $\alpha_i > 0$, except that the
framing of the meridian component is $-r_i$ instead of $r_i$.

\begin{figure}
\centerline{\includegraphics{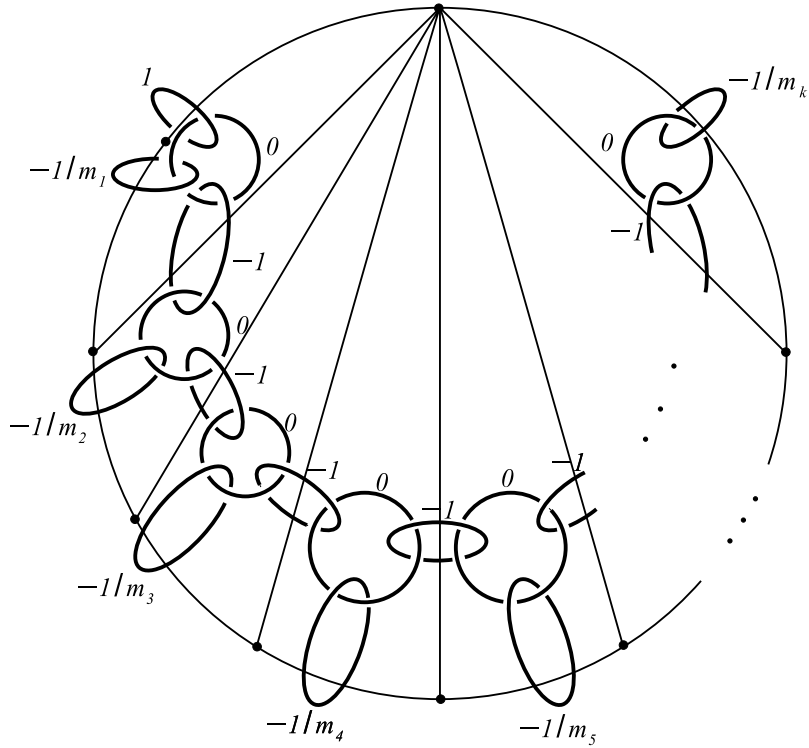}} \caption{The framed link
$K_i$ when $\alpha_i \le 0$.} \label{fig:glinks0}
\end{figure}

\begin{figure}
\centerline{\includegraphics{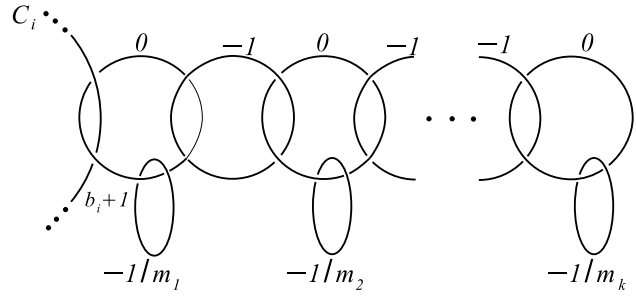}} \caption{Part of the
framed link $J$.} \label{fig:linksum}
\end{figure}

\begin{figure}
\centerline{\includegraphics{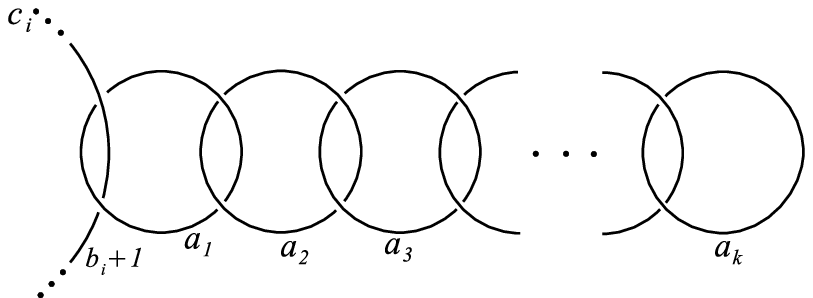}} \caption{Simplifying the
framed link $J$.} \label{fig:glinks1}
\end{figure}

\begin{figure}
\centerline{\includegraphics{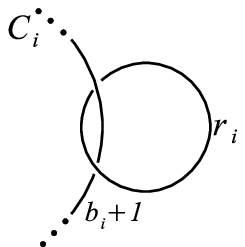}} \caption{Simplifying the
framed link $J$.} \label{fig:glinks2}
\end{figure}

Finally, one performs a slam dunk for each $i\in\cN$.  If $\za_i\le
0$, then the framing of $C_i$ becomes
  \begin{equation*}
b_i+\mul(E_i)^{-1}-\frac{1}{r_i}=b_i+1-(1-\za_i)=b_i+\za_i.
  \end{equation*}
If $\za_i>0$, then we have $b_i-1+(1+\za_i)=b_i+\za_i$.  So Dehn
surgery on the framed link $L$ is a bitwist manifold.

\end{proof}

\section{Realizing 3-manifolds as bitwist manifolds}
\label{sec:allarebitwist}\nosubsections

In this section we show that every closed connected orientable
3-manifold is a bitwist manifold.

Let $B$ be a braid with $n$ strands. Following \cite{PS}, we consider
the strands of $B$ as joining the points $A_i = (i,0,0)$ and $B_i =
(i,0,1)$ in $\bR^3$, $1\le i\le n$. The \textbf{closure} of $B$ is a link
in $S^3 = \bR^3 \cup \{\infty\}$ obtained by joining each $A_i$ and
$B_i$ by an arc such that the projections of these arcs on the
$xz$-plane are disjoint from each other and from the projection of $B$
onto the $xz$-plane. By a \textbf{generalized closure}, we only assume
that the endpoints $\{A_i\co 1\le i \le n\} \cup \{B_i\co 1\le i \le
n\}$ are joined by $n$ arcs whose projections are disjoint from each
other and from the projection of $B$.  This agrees with the definition
of closure given in \cite{Ka}, but is more restrictive than that
because we are not allowing any more crossings in the projection.

\begin{lemma}\label{lem:genclos}
Every link $L$ is a generalized closure of a pure braid.
\end{lemma}
\begin{proof}
Let $L$ be a link in $\bR^3$, and let $\pi\co\bR^3 \to \bR$ be the
projection onto the third coordinate. Then $L$ can be isotoped so
that, for some integer $n$, the height function on $L$ has $n$ local
maxima, which lie in $\pi^{-1}((1,\infty))$ and $n$ local minima,
which lie in $\pi^{-1}((-\infty,0))$. Furthermore, we can assume
that $L$ intersects the $xy$-plane in the points $A_i = (i,0,0)$,
$1\le i \le 2n$ , $L$ intersects the plane $z=1$ in the points $B_i
= (i,0,1)$, $1\le i \le 2n$, and all crossings of the projection of
$L$ onto the $xz$ -plane lie in $\pi^{-1}([0,1])$. (This follows,
for example, from Alexander's theorem, which states that $L$ can be
represented as the closure of an $n$-strand braid.) For convenience,
we call the components of $L \cap \pi^{-1}([0,1])$ the \textbf{strands}
of $L$, we call the components of $L \cap \pi^{-1}([1,\infty))$ the
\textbf{tops} of $L$, and we call the components of $L \cap
\pi^{-1}((-\infty,0])$ the \textbf{bottoms} of $L$. We first isotope
$L$ to a link $L_1$ so that there is a strand of $L_1$  joining
$A_1$ and $B_1$ and so that there is a top of $L_1$ joining $B_1$
and $B_2$. This can be done by sliding tops past each other and
possibly introducing a crossing in the projection of one top to
change the order of its endpoints in the projection. If the strand
of $L_1$ descending from $B_2$ ends at $A_2$, then we repeat this
process starting with the strand rising from $A_3$. Otherwise, by
sliding bottoms of $L_1$ past each other and possibly adding a
crossing in the projection of one bottom of $L_1$, we can isotope
$L_1$ to a link $L_2$ such that there is a strand of $L_2$ joining
$A_1$ and $B_1$, there is a top of $L_2$ joining $B_1$ and $B_2$,
there is a strand of $L_2$ joining $B_2$ and $A_2$, and there is a
bottom of $L_2$ joining $A_2$ and $A_3$. One next considers the
strand rising from $A_3$. One can continue this process to isotope
$L$ to a generalized closure of a pure braid with $2n$ strands.

\end{proof}

We next consider generators for the pure braid group. Let $K_n$ be
the pure braid group of isotopy classes of $n$-stranded pure braids.
Given $1\le i < j \le n$, let $b_{i,j}$ be the pure braid obtained
by doing a full twist on the collection of strands from the
$i^{\text{th}}$ to the $j^{\text{th}}$. Then (if the directions of
twisting are chosen properly) $a_{ij} = b_{i,j}b_{i+1,j}^{-1}$
is a pure braid for which the $i^{\text{th}}$ strand goes in front
of the $k^{\text{th}}$ strands, $i < k \le j$,  and then behind the
$k^{\text{th}}$ strands, $i < k \le j$. Since the elements $a_{ij}$,
$1\le i < j \le n$, generate the pure braid group, the elements
$b_{i,j}$, $1\le i < j \le n$, generate the pure braid group.

\begin{figure}
\centerline{\includegraphics{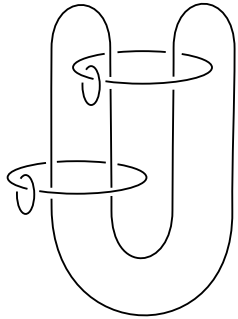}} \caption{The link
$L'$ for a simple example.} \label{fig:braidlink}
\end{figure}

\begin{figure}
\centerline{\includegraphics{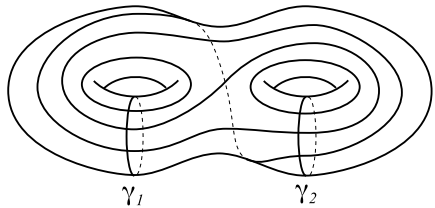}} \caption{The surface $S$
with meridians and nonmeridian link components.} \label{fig:curves}
\end{figure}

\begin{figure}
\centerline{\includegraphics{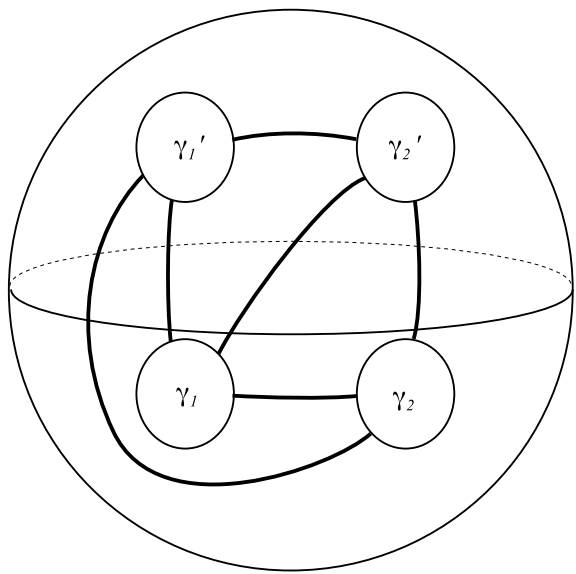}} \caption{Cutting open
the surface $S$ to get $S'$.} \label{fig:punctdisk}
\end{figure}

\begin{figure}
\centerline{\includegraphics{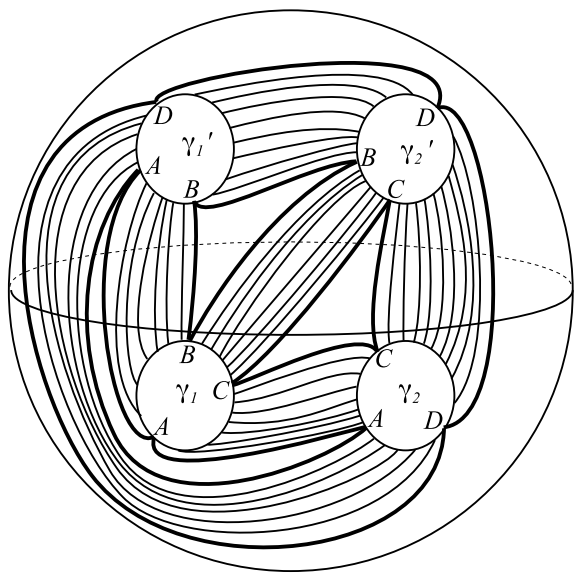}} \caption{Constructing the
faceted 3-ball.} \label{fig:fatarcs}
\end{figure}

\begin{figure}
\centerline{\includegraphics{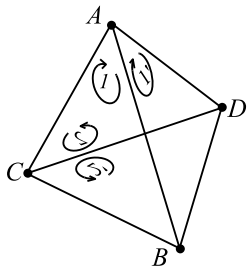}} \caption{The face-pairing
on the faceted 3-ball.} \label{fig:tet1p}
\end{figure}

\begin{thm}\label{thm:allm}
Every closed connected orientable 3-manifold is a bitwist
3-manifold.
\end{thm}
\begin{proof} Suppose $M$ is a closed connected orientable
3-manifold. By the Dehn-Lickorish theorem, $M$ can be obtained by Dehn
surgery on a framed link $L$. By Theorem~\ref{thm:consum} we can
assume that $L$ is not a split link.  By Lemma~\ref{lem:genclos}, $L$
is a generalized closure of a pure braid $B$.  We write $B$ in terms
of the generators $b_{i,j}$. We now view our projection of $L$ as
lying in the plane $\bR^2$. We view the braid $B$ as lying in a
rectangle $R$, with its strands joining the top and the bottom.  The
generators of $B$ lie in subrectangles which stack together to give
the rectangle $R$. Choose such a subrectangle $Q_g$ corresponding to a
generator $g$ of $B$. See Figure~\ref{fig:qgandrg}.a). Then $g$ is a
full twist on a set of consecutive strands of the braid in $Q_g$.  Let
$R_g$ be a subrectangle of $Q_g$ which contains only the consecutive
strands that are twisted in $g$. We next attach a rectangular block to
$R$ so that the bottom of the block is on $R_g$. The side of the block
facing the top of $R$ is the \textbf{front} of the block, and the side
of the block facing the bottom of $R$ is the \textbf{back} of the
block. We replace the strands of $R_g$ that are twisted by parallel
strands that go over the front of the block, along the top of the
block, and then back down the back of the block. We also drill out a
hole in the block that goes through the sides. See
Figure~\ref{fig:qgandrg}.b). In effect, we have added a handle to the
surface, and have replaced $g$ by a trivial braid which goes over the
handle. We also choose a circle for the boundary of the block's hole,
and we choose a meridian for the handle. We expand the meridian
slightly so that it links the arcs that go over the handle and the
circle in the boundary of the hole. See Figure~\ref{fig:qgandrg}.c).
We choose framing $0$ for the meridian, and framing $\pm 1$ (depending
on the direction of twist of the generator) for the circle in the
boundary of the hole. We shrink the block slightly so that blocks
corresponding to different subrectangles are disjoint.  Doing this for
each generator while maintaining the framings of the components of $L$
gives a framed link $L'$.  Let $S$ be the surface obtained from the
2-sphere $\bR^2 \cup \{\infty\}$ by adding a handle as described above
for each generator of $B$.

\begin{figure}
\centerline{\includegraphics{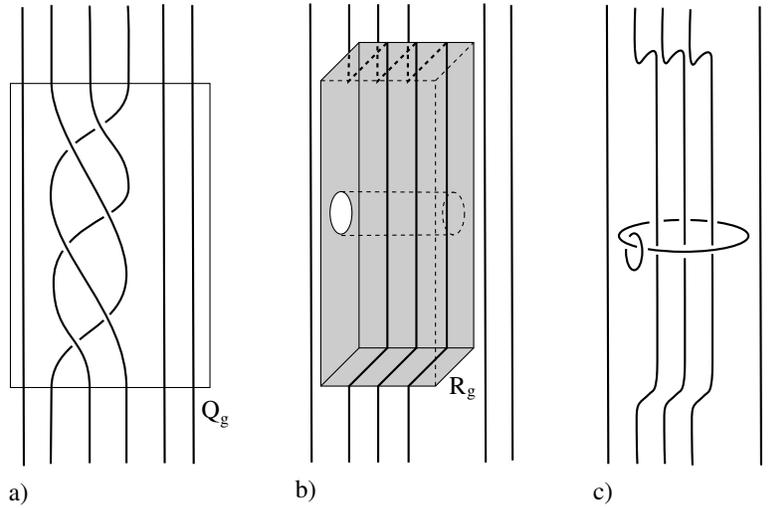}} \caption{Steps in the
construction of $L'$.} \label{fig:qgandrg}
\end{figure}

If we perform a slam-dunk on each circle along the boundary of a
hole, then the effect on $L'$ is to delete those circles and to
change the framing on each of the meridian circles to $\pm 1$. If we
now perform a Rolfsen twist along each of the meridian circles, then we
recover the original link $L$, but with framings changed by sums and
differences of squares of linking numbers of the meridian curves and
the components of $L$. Hence if we change the framings on $L'$  by
adding an appropriate integer to each of the components of $L$, we
get a framed link $L''$ such that $M$ is obtained from the 3-sphere
by surgery on $L''$. By Theorem~\ref{thm:glinks}, to prove
Theorem~\ref{thm:allm} it suffices to prove that $L'$ is a corridor
complex link whose face components are the meridians.

To get a face pairing, we cut open the surface $S$ along the
meridians. If there are $n$ meridians, the result is a 2-sphere
with $2n$ paired holes and disjoint arcs joining their boundaries.
We attach a disk to every hole to obtain a 2-sphere $S'$.  Since
$L$ is not a split link, the connected components of the
complement in $S'$ of the union of the arcs and closed disks are
all simply connected. The link in Figure~\ref{fig:braidlink} gives
rise to the surface with curves in Figure~\ref{fig:curves} (which
is taken from \cite{heegaard}).  Figure~\ref{fig:punctdisk} shows
the result $S'$ of cutting open $S$ and attaching disks.  We
fatten each arc to a quadrilateral, foliated by arcs parallel to
the core arc, so that adjacent quadrilaterals touch on the
boundaries of the $2n$ disks. See, for example,
Figure~\ref{fig:fatarcs}. We now collapse to a point each leaf in
a quadrilateral and the closure of each region in the complement
of the union of the paired disks and foliated quadrilaterals. By
Moore's theorem the quotient space $S''$ is a 2-sphere, with a
cell structure that consists of a vertex for each collapsed
complementary region, an edge for each collapsed foliated
quadrilateral, and a face for each of the $2n$ paired disks. We
define a face-pairing on the quotient space $S''$ in a
straightforward way. This defines a face-pairing $\ze$ for a
faceted 3-ball $P$ whose boundary is the 2-sphere $S''$. For the
example above, this is shown in Figure~\ref{fig:tet1p}. By
construction, $L'$ is a corridor complex link for $(P,\ze)$.

\end{proof}

\end{document}